\documentclass[11pt]{amsart}
\usepackage{amsmath,amssymb,mathrsfs,amscd}
\usepackage{graphicx}

\title{On a construction of the twistor spaces of Joyce metrics, II}
\author{Nobuhiro Honda
}  

\newcommand{\ol}{\overline}
\newcommand{\ra}{\rightarrow}
\newcommand{\lra}{\longrightarrow}

\newcommand{\Ra}{\Rightarrow}
\newcommand{\set}{\,|\,}
\newcommand{\proofend}{\hfill$\square$}

\newcommand{\vsp}{\vspace{3mm}}

\newtheorem{prop}{Proposition}[section]
\newtheorem{lemma}[prop]{Lemma}

\newtheorem{thm}[prop]{Theorem}
\newtheorem{rmk}[prop]{Remark}

\begin{document}

\maketitle

\begin{abstract} 
We explicitly construct the twistor spaces of Joyce metrics with torus action that are not treated  in Part I. 
This finishes a construction of all the twistor spaces of Joyce metrics on the connected sum of four complex projective planes.
\end{abstract}


\section{Introduction} 
This paper is a continuation of  \cite{H-1} (which we refer as Part I) on the twistor spaces of Joyce metrics on the connected sum of four complex projective planes.
In  Part I the twistor spaces  were explicitly constructed  for one family of Joyce metrics on $4\mathbf P^2$.
The purpose of the present paper is to do the same thing for another family of Joyce metrics on $4\mathbf P^2$.

The method   and ideas we will take is broadly similar to that of Part I.
Namely we first investigate the anticanonical system of the twistor space and show that the associated map is bimeromorphic onto its image.
Then the defining equations of the image are explicitly given in some $\mathbf P^6$-bundle over a projective line.
Finally we construct the biholomorphic model of the twistor space by taking  blowing-ups and blowing-downs of the anticanonical model.

However, there are many differences between the present case and Part I.
First, in Part I, the anticanonical model of the twistor space was a complete intersection in $\mathbf P^6$ of degree 8.
In the present case it is not a complete intersection  but an (ideal theoretic) intersection of 10 quadratic hypersurfaces in $\mathbf P^8$.
At first glance this seems to make  situation more complicated.
But it will turn out that  this is not true;
in Part I the anticanonical model  has `redundant divisors' which do not exist in the actual twistor space, and 
to blow-down these divisors are the most difficult step in the construction.
For the present case there are no such divisors and consequently the construction of the twistor spaces becomes considerably simpler.

We again remark that a reader who is interested in the construction itself can consult Section 3 and figures there  directly. 
The results of Section 2 are needed in order to prove that the constructed 3-folds are actually the twistor spaces of Joyce metrics.

\section{Defining equations of the anticanonical model}
As briefly explained in \S 2.1 of Part I, there are precisely 3 equivalence classes of smooth effective $K$-actions on $4\mathbf P^2$ and for each equivalence class  there associates a family of Joyce metrics on $4\mathbf P^2$.
One of the family consists of   special forms of LeBrun metrics \cite{LB91} and their twistor spaces are explicitly constructed in \cite{LB91}. The remaining 2 (equivalence classes of) smooth effective $K$-actions can be represented by the following sequences of isotropy subgroups respectively:
$$\{ K(1,0),K(1,1), K(0,1), K(-1,2), K(-2,3),K(-1,1)\}\,\cdots\text{Type I},$$
$$\{K(1,0),K(1,1),K(0,1),K(-1,2),K(-1,1), K(-2,1)\}\,\cdots\text{Type II}.$$
In  \cite{H-1} we constructed the twistor spaces of Joyce metrics on $4\mathbf P^2$ whose $K$-action is of Type I.
The purpose of the present paper is to do the same thing for Joyce metrics on $4\mathbf P^2$ whose $K$-action is of Type II.
Before its investigation, partly in order to set the notations, we again quote a result of A. Fujiki about structure of the twistor spaces of Joyce metrics in general:

\begin{prop}\label{prop-01}\cite{F00}
Let $Z$ be the twistor space of a Joyce metric on $n\mathbf P^2$ which is different from LeBrun metrics. Then  the following holds.
(i) $\dim |(-1/2)K_Z|=1$ and its members are  $G$-invariant.
Moreover, general members of the pencil are irreducible and smooth, and they become a toric surface with $K^2=8-2n$.
(ii) Let $S$ be a general member of the pencil 
and  $C$ the  unique $G$-invariant anticanonical curve on $S$. Then {\rm Bs} $|(-1/2)K_Z|=C$.
(iii) If we write $C=\sum_{i=1}^{2n+4}C_i$ in such a way that $C_i$ and $C_{i+1}$ intersect, the real structure on $Z$ exchanges $C_i$ and $C_{i+n+2}$, where the subscriptions are counted modulo $2n+4$.
(iv) The pencil $|(-1/2)K_Z|$ has precisely $(n+2)$ singular members and all of them consist of two smooth irreducible toric surfaces that are conjugate of each other.
(v) If we write $S_i=S_i^++S_i^-$ ($1\le i\le n+2$) for the reducible members, $S^+_i$ and $S^-_i$ divide $C$ into 'halves' in the sense that both $S_i^+\cap C$ and $S_i^-\cap C$ are connected (cf.\! (iii)). 
Moreover,  $L_i:=S^+_i\cap S_i^-$ is a $G$-invariant twistor line. \end{prop}

\begin{figure}
\includegraphics{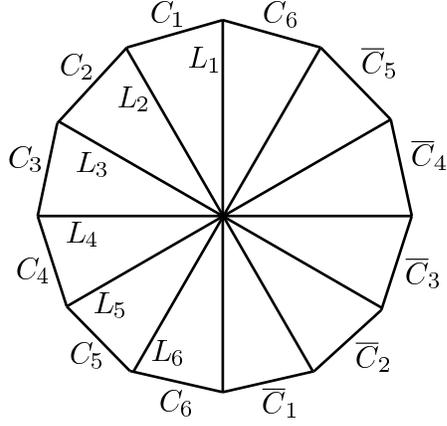}
\caption{The numbering for $G$-invariant twistor lines}
\label{fig-6tls}
\end{figure}

In the following $Z$ always denote the twistor space of a Joyce metric on $4\mathbf P^2$ whose $K$-action is Type II.
In order to investigate the anticanonical system of $Z$, as in the case of Type I,  we study bi-anticanonical system of a general member $S$ of the pencil $|(-1/2)K_Z|$.
By Proposition \ref{prop-01}, $S$  is a smooth toric surface with $K^2=0$.
For the cycle of $G$-invariant anticanonical curve $C$ of $S$, we again write in the form
$$
C=\sum_{i=1}^6C_i+\sum_{i=1}^6\ol{C}_i,
$$
where $\ol{C}_i$ denotes the conjugation of $C_i$ and the components satisfy
 $C_iC_{i+1}= \ol{C}_i\ol{C}_{i+1}=1$ for $1\le i\le 5$ and $C_6\ol{C}_1=\ol{C}_6 C_1=1$ (cf. Figure \ref{fig-6tls}).
Moreover, after possible renumbering, the self-intersection numbers of $C_i$ inside $S$ can be supposed to satisfy, this time
$$
C_1^2=-3,\,C_2^2=-1,\,C_3^2=-3,\,C_4^2=-1,\,C_5^2=-3,\,C_6^2=-1.
$$


Then we have the following result on the bi-anticanonical system of $S$:

\begin{lemma}\label{lemma-01}
Let $S$ be an irreducible member of the pencil $|(-1/2)K_Z|$  and $C=\sum C_i+\sum \ol{C}_i$ the  unique $G$-invariant anticanonical curve on $S$ as above.
Then we have the following.
(i) {\rm{Bs}}\,$|-2K_S|=C_1+C_3+C_5+\ol{C}_1+\ol{C}_3+\ol{C}_5$ and the movable part of the system  is base point free and 6-dimensional.
(ii) Its associated morphism $\phi:S\ra\mathbf P^6$ is birational onto its image.
(iii) The image of $\phi$ is defined, ideal theoretically, by the following 9 quadratic equations
\begin{align}
x_1x_2&=ax_0^2 &ax_3x_5&=dx_0x_2 & dx_4x_6&=bcx_0x_1\notag\\ 
x_3x_4&=bx_0^2& bx_1x_5&=dx_0x_4 &dx_2x_6&=acx_0x_3\label{toric01}\\
x_5x_6&=cx_0^2&cx_1x_3&=dx_0x_6& dx_2x_4&=abx_0x_5 \notag
\end{align}
where  $(x_0,\cdots,x_6)$ is a homogeneous coordinate on $\mathbf{P}^6$, and $a,b,c,d$ are non-zero constants.
(iv) If $S$ is a real member of the pencil, the real structure on \,$\mathbf P^6=\mathbf P^{\vee}H^0(-2K_S)$ induced from $S$ is given by, in the above coordinate, 
\begin{equation}\label{rs01}
(x_0,x_1,x_2,x_3,x_4,x_5,x_6)\longmapsto
(\ol{x}_0,\ol{x}_2,\ol{x}_1,\ol{x}_4,\ol{x}_3,\ol{x}_6,\ol{x}_5).
\end{equation}
Moreover, we have $a,b,c\in\mathbf R$ and $\,abc=|d|^2$.
\end{lemma}

\noindent 
Proof. The following argument is similar to that of Lemma 2.3 in \cite{H-1}.
Actually (i) and (ii) can be verified by  standard calculations as in \cite{H-1}.
For (iii), we again use the theory of toric varieties to compute the natural $G$-action on $H^0(-2K_S)^{\vee}\simeq\mathbf C^7$ and conclude that it is given by
\begin{equation}
\label{act01}
(x_0,x_1,x_2,x_3,x_4,x_5,x_6)\longmapsto (x_0, sx_1,s^{-1}x_2,tx_3,t^{-1}x_4,s^{-1}t^{-1}x_5,stx_6)\hspace{2mm}\text{for}\hspace{2mm}(s,t)\in G
\end{equation}
for a coordinate $(x_0,x_1,\cdots,x_6)$ on  $H^0(-2K_S)^{\vee}$.
By similar computation we obtain \eqref{rs01} of (iv).
To obtain defining equations of the image $\phi(S)$ in $\mathbf P^6$, we consider a rational map $\mathbf P^6\ra\mathbf P^3\times\mathbf P^1$ defined by 
\begin{equation}\label{quotient01}
(x_0,x_1,x_2,x_3,x_4,x_5,x_6)\longmapsto (x_0^2,x_1x_2, x_3x_4,x_5x_6)\times (x_0^3, x_1x_3x_5).
\end{equation}
It is easily seen that the map \eqref{quotient01} is surjective and $G$-equivariant, where the $G$-action on $\mathbf P^3\times\mathbf P^1$  is the trivial action.
On the other hand, the morphism $\phi:S\ra\mathbf P^6$ is $G$-equivariant.
Hence $\phi(S)$ must be contained in some fiber of the rational map \eqref{quotient01}.
Let $(y_0,y_1,y_2,y_3)\times(z_0,z_1)\in\mathbf P^3\times\mathbf P^1$  be  the point of which the fiber contains $\phi(S)$.
By \eqref{quotient01} if $y_0=0$, then $\phi(S)$ would be contained in a hyperplane $x_0=0$.
This cannot happen since the image must be non-degenerate.
By the same reason, we have $y_i\neq 0$ and $z_j\neq 0$ for all $i=1,2,3$ and $j=0,1$.
So we can suppose that $(y_0,y_1,y_2,y_3)=(1,a,b,c)$ and $(z_0,z_1)=(1,d)$ for some $a,b,c,d\in\mathbf C^*$.
Namely the image $\phi(S)$ is contained in 
\begin{equation}\label{defeq01}
x_1x_2=ax_0^2,\,\,x_3x_4=bx_0^2,\,\,x_5x_6=cx_0^2,\,\,x_1x_3x_5=dx_0^3.
\end{equation}
\eqref{defeq01} defines a surface in $\mathbf P^6$ but it is not irreducible.
Actually, it is immediate to see that the surface \eqref{defeq01} contains 7 planes (like $x_0=x_1=x_3=x_5=0$).
We can deduce by elementary computations that there is a unique non-degenerate  component and  its defining ideal is generated by 9 equations \eqref{toric01}. 
The final assertion of (iv) can be verified by direct calculation using \eqref{rs01}.
\proofend

\vsp
We note that the degree of the surface \eqref{toric01} is six, since $\phi$ is a birational morphism and the self-intersection number of the movable part of $|-2K_S|$ is six.
Namely, the bi-anticanonical model of our toric surface $S$ is a sextic surface  defined by the 9 quadratic equations \eqref{toric01}.

\begin{lemma}\label{lemma-011} 
(i) If $abcd\neq 0$, then the 9 quadratic equations  \eqref{toric01} define a smooth toric surface in $\mathbf P^6$ satisfying $K^2=6$.
(ii) If $a,b,c,d$ satisfy $(a,b,c)\in\mathbf R^3$ and $|d|^2=abc$, then the subvariety in $\mathbf P^6$ defined by  \eqref{toric01}, which becomes automatically a real subset in $\mathbf P^6$ with respect to the real structure  \eqref{rs01}, has a real point iff $a\ge0, b\ge 0$ and $c\ge 0$ hold.
(iii) If  $\phi:S\ra\mathbf P^6$ is as in Lemma \ref{lemma-01}, then $\phi:S\ra\phi(S)$ is precisely the blowing-downs of the six $(-1)$  curves $C_2,C_4, C_6,\ol{C}_2,\ol{C}_4$ and $\ol{C}_6$.
\end{lemma}

\noindent
The proof is not difficult and we leave it to the reader.

\begin{lemma}\label{lemma-02}
Let $Z$ be the twistor space of a Joyce metric on $4\mathbf P^2$ whose $K$-action is Type II.
Then we have $\dim H^0(-K_Z)=9$. 
Further, {\rm{Bs}}\,$|-K_Z|=C_1+C_3+C_5+\ol{C}_1+\ol{C}_3+\ol{C}_5$.
\end{lemma}

\noindent
Proof.
By the same reason as in the proof of Theorem 2.2 in \cite{H-1}, we have a $K$-equivariant exact sequence 
\begin{equation}\label{ses01}
0\lra H^0((-1/2)K_Z)\lra H^0(-K_Z)\lra H^0(-2K_S)\lra 0.
\end{equation}
We have $\dim H^0((-1/2)K_Z)=2$, and $\dim H^0(-2K_S)=7$ by Lemma \ref{lemma-01}.
Hence we obtain $\dim H^0(-K_Z)=9$.
By the surjectivity of the restriction map $H^0(-K_Z)\lra H^0(-2K_S)$ in \eqref{ses01}, we obtain Bs\, $|-K_Z|=$\,Bs\,$|-2K_S|$.
The latter is $C_1+C_3+C_5+\ol{C}_1+\ol{C}_3+\ol{C}_5$ by Lemma \ref{lemma-01} (i).
\proofend

\vsp
The surjectivity of the restriction map $H^0(-K_Z)\lra H^0(-2K_S)$ in \eqref{ses01} again becomes crucial in our description of the anticanonical model of $Z$.
By the same reasoning as in the proof of \cite[Theorem 2.2 (ii)]{H-1}, the $G$-action naturally induced on $H^0((-1/2)K_Z)\simeq\mathbf C^2$ is trivial.
Then the exact sequence \eqref{ses01} and the action \eqref{act01} on $H^0(-2K_S)^{\vee}$ mean that the natural $G$-action on $H^0(-K_Z)^{\vee}\simeq\mathbf C^9$ is given by
\begin{equation}\label{act02}
(x_1,\cdots,x_9)\mapsto (sx_1,s^{-1}x_2,tx_3,t^{-1}x_4,s^{-1}t^{-1}x_5,stx_6,x_7,x_8,x_9)
\end{equation}
for a coordinate $(x_1,\cdots,x_9)$ on $H^0(-K_Z)^{\vee}$.
As in  the proof of  \cite[Theorem 2.2 (ii)]{H-1}, let $V$ be the 3-dimensional vector subspace of $H^0(-K_Z)$ generated by the image of the natural map
\begin{equation}\label{f01}
H^0((-1/2)K_Z)\times H^0((-1/2)K_Z)\lra H^0(-K_Z);\hspace{3mm}
(\sigma_1,\sigma_2)\longmapsto \sigma_1\otimes\sigma_2.
\end{equation}
Then $V$ coincides with  the vector subspace of $H^0(-K_Z)$ consisting of $G$-fixed elements.
If $\Phi:Z\ra\mathbf P^{\vee}H^0(-K_Z)=\mathbf P^8$ and $\Psi:Z\ra \Lambda:=\mathbf P^{\vee}H^0((-1/2)K_Z)=\mathbf P^1$ denote the meromorphic maps associated to $|-K_Z|$ and $|(-1/2)K_Z|$ respectively, we obtain a commutative diagram of meromorphic maps
\begin{equation}\label{cd01}
 \CD
 Z@>\Phi>> \mathbf P^8\\
 @V\Psi VV @VV\pi V\\
 \Lambda@>\iota>> \mathbf P^2
 \endCD
 \end{equation}
 where $\mathbf P^2=\mathbf P^{\vee}V$ and
$\iota:\Lambda\ra \mathbf P^2=\mathbf P^{\vee}V$ is an inclusion as a conic, as in  \cite[the diagram (11)]{H-1}, and $\pi:\mathbf P^8\ra\mathbf P^2$ is the projection determined by the inclusion $V\subset H^0(-K_Z)$. 
For $\lambda\in\Lambda$, we denote $S_{\lambda}:=\Psi^{-1}(\lambda)\in |(-1/2)K_Z|$.
By the surjectivity of the restriction map in the sequence \eqref{ses01}, the restriction $\Phi|_{S_{\lambda}}$ is precisely the rational map $\phi_{\lambda}:S_{\lambda}\ra \mathbf P^{\vee}H^0(-2K_{S_{\lambda}})=\pi^{-1}(\iota(\lambda))$ associated to $|-2K_{S_{\lambda}}|$.
Therefore $\Phi$ is bimeromorphic onto its image since $\iota$ in \eqref{cd01} is an embedding.
Let $\mathbf P^5_{\infty}$ be the center of the projection $\pi$.
We take the blowing-up of $\mathbf P^8$ with center $\mathbf P^5_{\infty}$.
The resulting space is the total space of a $\mathbf P^6$-bundle $\mathbf P(\mathscr O(1)^{\oplus 6}\oplus \mathscr O)\ra\mathbf P^2$.
Restricting onto the conic $\Lambda\subset\mathbf P^2$, we obtain a $\mathbf P^6$-bundle $\mathbf P(\mathscr O(2)^{\oplus 6}\oplus \mathscr O)\ra\Lambda=\mathbf P^1$.
Then the meromorphic map $\Phi$ naturally lifts as a meromorphic map $$\Phi_1:Z\ra\mathbf  P(\mathscr O(2)^{\oplus 6}\oplus \mathscr O).$$
$\Phi_1$ is bimeromorphic onto its image since $\Phi$ is.

In the following we determine defining equations of the image $\Phi_1(Z)$ in the total space of 
$\mathbf P(\mathscr O(2)^{\oplus 6}\oplus \mathscr O)$. 
To write down the equations explicitly, we first give a numbering for the set of $G$-invariant twistor lines, as in Figure \ref{fig-6tls}.
Then we obtain a  numbering for the 6 real points in $\Lambda$.
Let $\lambda_1,\cdots,\lambda_6$ be the real points (real numbers) obtained in this way.
Obviously either $\lambda_1<\cdots<\lambda_6$ or $\lambda_1>\cdots>\lambda_6$ holds.

\begin{thm}\label{thm-01}
Let $Z$ be the twistor space of a Joyce metric on $4\mathbf P^2$ whose $K$-action is Type II and   $\{\lambda_1,\lambda_2\cdots,\lambda_6\}$ real numbers obtained as above.
Let $\Phi:Z\ra \mathbf P^8$ be the meromorphic map associated to the anticanonical system of \,$Z$ and $\Phi_1:Z\ra\mathbf P(\mathscr O(2)^{\oplus 6}\oplus \mathscr O)$  the natural lift of $Z$ as above. 
Then $\Phi_1$ is bimeromorphic onto its image, and the image $\Phi_1(Z)$ is defined by the following 9 equations
\begin{align}
x_1x_2&=-(\lambda-\lambda_2)(\lambda-\lambda_5)(\lambda-\lambda_3)(\lambda-\lambda_6)\,x_0^2,\label{defeq02-1}\\
x_3x_4&=-(\lambda-\lambda_1)(\lambda-\lambda_4)(\lambda-\lambda_2)(\lambda-\lambda_5)\,x_0^2,\label{defeq02-2}\\
x_5x_6&=\hspace{3mm}(\lambda-\lambda_1)(\lambda-\lambda_4)(\lambda-\lambda_3)(\lambda-\lambda_6)\,x_0^2,\label{defeq02-3}\\
x_3x_5&=-(\lambda-\lambda_1)(\lambda-\lambda_4)\,x_0x_2,\label{defeq02-4}\\
x_5x_1&=-(\lambda-\lambda_3)(\lambda-\lambda_6)\,x_0x_4,\label{defeq02-5}\\
x_1x_3&=\hspace{3mm}(\lambda-\lambda_2)(\lambda-\lambda_5)\,x_0x_6,\label{defeq02-6}\\
x_4x_6&=-(\lambda-\lambda_1)(\lambda-\lambda_4)\,x_0x_1\label{defeq02-7}\\
x_2x_6&=-(\lambda-\lambda_3)(\lambda-\lambda_6)\,x_0x_3,\label{defeq02-8}\\
x_2x_4&=\hspace{3mm}(\lambda-\lambda_2)(\lambda-\lambda_5)\,x_0x_5,\label{defeq02-9}
\end{align}
where $\lambda$ is an affine coordinate on $\Lambda$, and  $(x_1,x_2,\cdots,x_6;x_0)$ is a fiber coordinate of the bundle $\mathscr O(2)^{\oplus 6}\oplus \mathscr O\ra\Lambda$ which is valid on $\lambda\neq\infty$.
(Here, all the equations \eqref{defeq02-1}--\eqref{defeq02-9} take values in the bundle $\mathscr O(4)\ra\Lambda$.)
\end{thm}

\begin{rmk}{\em 
In \cite[Theorems 2.2 and 2.5]{H-1} we gave  equations of the anticanonical model of $Z$ (of Type I) defined in a projective space $\mathbf P^6$. 
However, in the present case, we do not give defining equations in $\mathbf P^8$  to avoid unnecessary complexity.
Actually, the above realization in the projective bundle is more natural and useful, and the explicit construction we will give later starts from the above model.
(This was already true in \cite{H-1}.)
}
\end{rmk}

For the proof of Theorem \ref{thm-01}, the most important  point is again to determine the image of reducible members of the pencil $|(-1/2)K_Z|$ under the meromorphic map $\Phi$.
To this end we first have the following lemma. 
For the natural numbering for the reducible members of  the pencil $|(-1/2)K_Z|$, we follow  Proposition \ref{prop-01} (v). 
Further, to make a distinction between the 2 components $S_i^+$ and $S_i^-$ by declaring that $S_i^+$ always contains $C_6$.

\begin{lemma}
\label{lemma-red}
The set of all $G$-invariant members of the anticanonical system $|-K_Z|$ consists of the following divisors;
(a) $S+S'$, where both $S$ and $S'$ are members of $|(-1/2)K_Z|$,
(b) $S_1^++S_3^++S_4^++S_6^+$, $S_1^++S_2^++S_4^++S_5^-$, $S_2^++S_3^-+S_5^-+S_6^-$, and the conjugations of these 3 divisors.
\end{lemma}

\noindent
Proof. 
By calculating the Chern classes of the divisors in the item $(b)$, we can verify that all the 6 divisors in $(b)$ belong to the anticanonical system of $Z$.
Thus to finish a proof it suffices to show that there are no $G$-invariant anticanonical divisor which are not of the form of $(a)$ and $(b)$.
But this fact immediately follows from the fact that 
the naturally induced $G$-action on $H^0(-K_Z)^{\vee}$ is given by \eqref{act02}.
\proofend

\vsp
In contrast to the twistor spaces studied in \cite{H-1}, the roles of all reducible members $S_i^++S_i^-$, $1\le i\le 6$, are basically equal:

\begin{lemma}\label{lemma-rest01}
Let $W_i^{\pm}$ $(1\le i\le 6$) be the vector space which is the image of the restriction map $H^0(-K_Z)\ra H^0(-K_Z|_{S_i^{\pm}})$.
Then we have the following.
(i) $\dim W_i^+=\dim W_i^-=5$ for any $1\le i\le 6$.
(ii) The linear systems associated to the vector spaces $W_i^+$ and $W_i^-$ have a unique base point.
(iii) The rational map associated to the linear system in (ii) is a birational map from $S_i^{\pm}$ onto a smooth cubic surface in $\mathbf P^4$, and it is biholomorphic on the twistor line $L_i$ in $S_i^{\pm}$.
\end{lemma}

\noindent
Proof.
Since the situation is basically the same for any $W_i^{\pm}$, we only consider $W_1^+$. 
To prove (i) it suffices to show $\dim H^0(-K_Z-S_1^+)=4$.
To calculate this cohomology group, we find $G$-invariant members of the system $|-K_Z-S_1^+|$.
By Lemma \ref{lemma-red} we obtain that the two divisors
$
S_3^++S_4^++S_6^+,\,\,S_2^++S_4^++S_5^-
$
are $G$-invariant members of  $|-K_Z-S_1^+|$.
Further, if $S$ is a member of $|(-1/2)K_Z|$, $S+S_1^-$ is also a member of $|-K_Z-S_1^+|$.
By verifying linear independence of these $G$-invariant members, we obtain that $\dim |-K_Z-S_1^+|\ge 3$.
But $\dim |-K_Z-S_1^+|=3$ must hold since  all $G$-invariant divisors in $Z$ are irreducible components of a member of $|(-1/2)K_Z|$ and hence the divisors of the above forms are all  $G$-invariant members of  $|-K_Z-S_1^+|$.
Therefore we obtain $\dim H^0(-K_Z-S_1^+)=4$ and (i) is proved.
By similar argument, we explicitly obtain all $G$-invariant members of $|-K_Z|$ which do not have $S_1^+$ as an irreducible component.
These members yield $G$-invariant members of  the system $|W_1^+|$ by restriction.
Consequently we concretely obtain all $G$-invariant members of $|W_1^+|$.
Once these members are obtained it is easy to determine the property of the rational maps associated to the linear system, and we conclude that (ii) and (iii) of the lemma hold.
We omit the detail.
\proofend

\vsp
\noindent
Proof of Theorem \ref{thm-01}.
The bimeromorphicity of $\Phi_1$ was already obtained.
If $\lambda\in\Lambda$ is a point such that $S_{\lambda}\in|(-1/2)K_Z|$ is irreducible (namely non-singular), then the image $\Phi_1(S_{\lambda})$ is defined by the 9 quadratic equations \eqref{toric01} in $\mathbf P^6$ as in Lemma \ref{lemma-01}.
Hence $\Phi_1(Z)$ must be of the form  \eqref{toric01}, where this time $a,b$ and $c$ are    quartic polynomials  of $\lambda$ and $d$ is a polynomial of degree 6 of $\lambda$, since $\Phi_1(Z)$ is in the bundle $\mathbf P(\mathscr O(2)^{\oplus 6}\oplus \mathscr O)\ra\mathbf P^1$.
We know by Lemma \ref{lemma-01} (iii) that if  $\lambda\neq\lambda_i$ for $1\le i\le 6$, then  $a(\lambda)b(\lambda)c(\lambda)d(\lambda)\neq 0$.
Moreover, taking the real structure into account, it follows by Lemma \ref{lemma-01} (iv) that 
\begin{equation}\label{reality01}
a(\lambda)b(\lambda)c(\lambda)=|d(\lambda)|^2\,\,\,{\rm{for}}\,\,\,\lambda\in\mathbf R.
\end{equation}
Hence if $d(\lambda_i)\neq 0$, then $a(\lambda_i)\neq 0$, $b(\lambda_i)\neq 0$ and $c(\lambda_i)\neq 0$.
By Lemma  \ref{lemma-011} this implies that  if $d(\lambda_i)\neq 0$, then $\Phi(S_{\lambda_i})$ would become a irreducible  sextic surface.
On the other hand, by  Lemma \ref{lemma-rest01},  $\Phi(S_i^+)$ and $\Phi(S_i^-)$ are cubic surfaces in $\mathbf P^4$.
This is a contradiction and hence we obtain that $d(\lambda_i)=0$ for any $1\le i\le 6$.
Thus since $\deg d(\lambda)=6$ we obtain 
$$
d(\lambda)=d_0\prod_{i=1}^6(\lambda-\lambda_i),
$$
for some constant $d_0$. Hence by \eqref{reality01} we obtain
\begin{equation}\label{reality02}
a(\lambda)b(\lambda)c(\lambda)=|d_0|^2\prod_{i=1}^6(\lambda-\lambda_i)^2.
\end{equation}
The problem is how to distribute the zeros $\{\lambda_i\set 1\le i\le 6\}$ to $a(\lambda)$, $b(\lambda)$ and $c(\lambda)$.
We first show that $a, b$ and $c$ do not have a double root.
Suppose for instance  that $\lambda_1$ is a double root of  $a$.
Then by \eqref{reality02} we have $b(\lambda_1)c(\lambda_1)\neq 0$.
Hence by the equations \eqref{toric01} the image $\Phi_1(S_1)$ is contained in the subvariety (in $\mathbf P^6$) 
\begin{equation}\label{toric02}
x_1x_2= x_1x_3=x_2x_6=x_0x_2=x_2x_4=x_1x_5=x_0x_1=0,\,\,
x_3x_4=bx_0^2, \,\,x_5x_6=cx_0^2
\end{equation}
where $b=b(\lambda_1)\neq 0$ and $c=c(\lambda_1)\neq 0$.
It can be readily seen that  \eqref{toric02} consists of a conjugate pair of 
 planes $\{x_0=x_1=x_4=x_6=0\}$, $\{x_0=x_2=x_3=x_5=0\}$ and a real quartic surface $\{x_1=x_2=0, x_3x_4=bx_0^2,x_5x_6=cx_0^2\}$.
On the other hand, we know by Lemma \ref{lemma-rest01} that $\Phi(S_1^+)$ and $\Phi(S_1^-)$ are cubic surfaces in $\mathbf P^6$.
This is a contradiction and hence we obtain that $a(\lambda)$ has no double root at $\lambda=\lambda_1$. 
Similar argument shows that $a(\lambda), b(\lambda)$ and $c(\lambda)$ has no double root at $\lambda=\lambda_i$ for any $1\le i\le 6$.

Next for simplicity of notations we set $A=\{\lambda_i\set a(\lambda_i)=0\}$, $B=\{\lambda_i\set b(\lambda_i)=0\}$ and $C=\{\lambda_i\set c(\lambda_i)=0\}$.
Of course we have $A\cup B\cup C=\{\lambda_i\set 1\le i\le 6\}$  as sets by \eqref{reality02}.
The fact that $a,b,c$ do not have a double root implies that $A\cap B$, $B\cap C$ and $C\cap A$ consist of precisely 2 elements.
Therefore we obtain an equality 
\begin{equation}\label{dcp01}
\{\lambda_1,\lambda_2,\lambda_3,\lambda_4,\lambda_5,\lambda_6\}=(A\cap B)\cup(B\cap C)\cup(C\cap A).
\end{equation}
From this, we can verify by direct calculations that all the six varieties in $\mathbf P^6=\Phi^{-1}(\lambda_i)$ defined by \eqref{toric01} consist of 2 irreducible components, and that both of the components are cubic surfaces contained in some linear $\mathbf P^4$ in $\mathbf P^6$.
Moreover, the intersection of the 2 cubic surfaces is an irreducible conic (contained in some $\mathbf P^2\subset\mathbf P^6)$.
 It can be derived by direct computations using the equations \eqref{toric01} and the $G$-action \eqref{act01} that if $A\cap B=\{\lambda_i,\lambda_j\}$ $(i\neq j)$, then the two conics over $\lambda_i$ and $\lambda_j$ have the same isotropy subgroup in $G$.
 On the other hand by Lemma \ref{lemma-rest01} the conic over $\lambda_i$ (which is the intersection of the 2 cubic surfaces) is a biholomorphic image of $L_i$ , where $L_i=S_i^+\cap S_i^-$ as before.
Moreover, we can readily calculate the isotropy subgroup along $L_i$ for $1\le i\le 6$ and deduce that $L_1$ and $L_4$ have the same  isotropy subgroup, and the same is true for $L_2$ and $L_5$, as well as $L_3$ and $L_6$.
These imply that the above decomposition \eqref{dcp01} must be  $\{\lambda_i\set 1\le i\le 6\}=\{\lambda_1,\lambda_4\}\cup\{\lambda_2,\lambda_5\}\cup\{\lambda_3,\lambda_6\}$.
Then by taking a cyclic permutation if necessary, we can suppose that $A=\{\lambda_2,\lambda_5,\lambda_3,\lambda_6\}$, $B=\{\lambda_1,\lambda_4,\lambda_2,\lambda_5\}$ and $C=\{\lambda_1,\lambda_4,\lambda_3,\lambda_6\}$.
Thus we obtain
$$a(\lambda)=a_0(\lambda-\lambda_2)(\lambda-\lambda_5)(\lambda-\lambda_3)(\lambda-\lambda_6)$$
$$b(\lambda)=b_0(\lambda-\lambda_1)(\lambda-\lambda_4)(\lambda-\lambda_2)(\lambda-\lambda_5)$$
$$c(\lambda)=c_0(\lambda-\lambda_1)(\lambda-\lambda_4)(\lambda-\lambda_3)(\lambda-\lambda_6).$$
for some $a_0,b_0,c_0\in\mathbf R^{\times}$.
Because $\Phi(S_{\lambda})$ has no real points for $\lambda\in\mathbf R$, Lemma \ref{lemma-011}  (ii) means that $a_0<0, b_0<0$ and $c_0>0$ hold.
Now we set
$$
\alpha_0=1,\,\,
\alpha_1=\sqrt{-a_0}\,e^{i\theta_1}x_1,\,\,
\alpha_3=\sqrt{-b_0}\,e^{i\theta_3}x_3,\,\,
\alpha_5=\sqrt{c_0}\,e^{i\theta_5},\,\,
\alpha_2=\ol{\alpha}_1,\,\,
\alpha_4=\ol{\alpha}_3,\,\,
\alpha_6=\ol{\alpha}_5
$$
where $\theta_1,\theta_3$ and $\theta_5$ are real numbers satisfying $d_0=\sqrt{a_0b_0c_0}\,e^{i(\theta_1+\theta_3+\theta_5)}$.
 (This is possible by the reality condition $|d_0|^2=a_0b_0c_0$).
 Next we take a  coordinate change  $X_i=\alpha_ix_i$, $0\le i\le 6$.
  Then dividing the both hand sides of the equations in \eqref{toric01} by the greatest common polynomials in $\lambda$,  we obtain that defining equations of the image $\Phi_1(Z)$ can be taken to be a normal form \eqref{defeq02-1}--\eqref{defeq02-9} as in the theorem.
  \proofend
  
\section{An explicit construction of the twistor spaces}
In this section we explicitly construct the twistor spaces of Joyce metrics whose $K$-action is of Type II.
As in the previous section let $\Lambda=\mathbf P^{\vee}H^0((-1/2)K_Z)\simeq\mathbf P^1$ be the projective line parametrizing the members of the pencil $|(-1/2)K_Z|$, and we write $S_{\lambda}\in |(-1/2)K_Z|$ for $\lambda\in\Lambda$. 
Let $\lambda_1,\cdots,\lambda_6\in\Lambda$ be points for which $S_{\lambda_i}=S_i\,(=S_i^++S_i^-)$ hold.
As explained in the first few paragraphs in  \cite[Section 3]{H-1}, by a result of Fujiki \cite[Theorem 9.1]{F00}, ${\lambda_1,\cdots,\lambda_6}$ can be supposed to coincide with the 6 real numbers involved in the construction of Joyce metrics which were written $p_1,\cdots,p_6$ in \cite[Theorem 3.3.1]{J95}.
(These numbers are called `conformal invariant' of Joyce metrics in \cite{F00}.)

Let $Z$ be the twistor space of a Joyce metric on $4\mathbf{CP}^2$ whose $K$-action is Type II and whose conformal invariant is $\lambda_1<\cdots<\lambda_6$.
By Theorem \ref{thm-01}, $Z$ is bimeromorphic to a 3-fold $X_1$ in the bundle $\mathbf P(\mathscr O(2)^{\oplus 6}\oplus \mathscr O)\ra\Lambda=\mathbf P^1$ defined by the equations \eqref{defeq02-1}--\eqref{defeq02-9}.
Let $p_1:X_1\ra\Lambda$ be the projection.
Note that there is a trivial subbundle 
$$
\mathbf P(\mathscr O(2)^{\oplus 6})\subset\mathbf P(\mathscr O(2)^{\oplus 6}\mathscr O)
$$
and the intersection of $X_1$ with this trivial subbundle is a product of $\mathbf P^1$ with a cycle of 6 rational curves.

\begin{figure}
\includegraphics{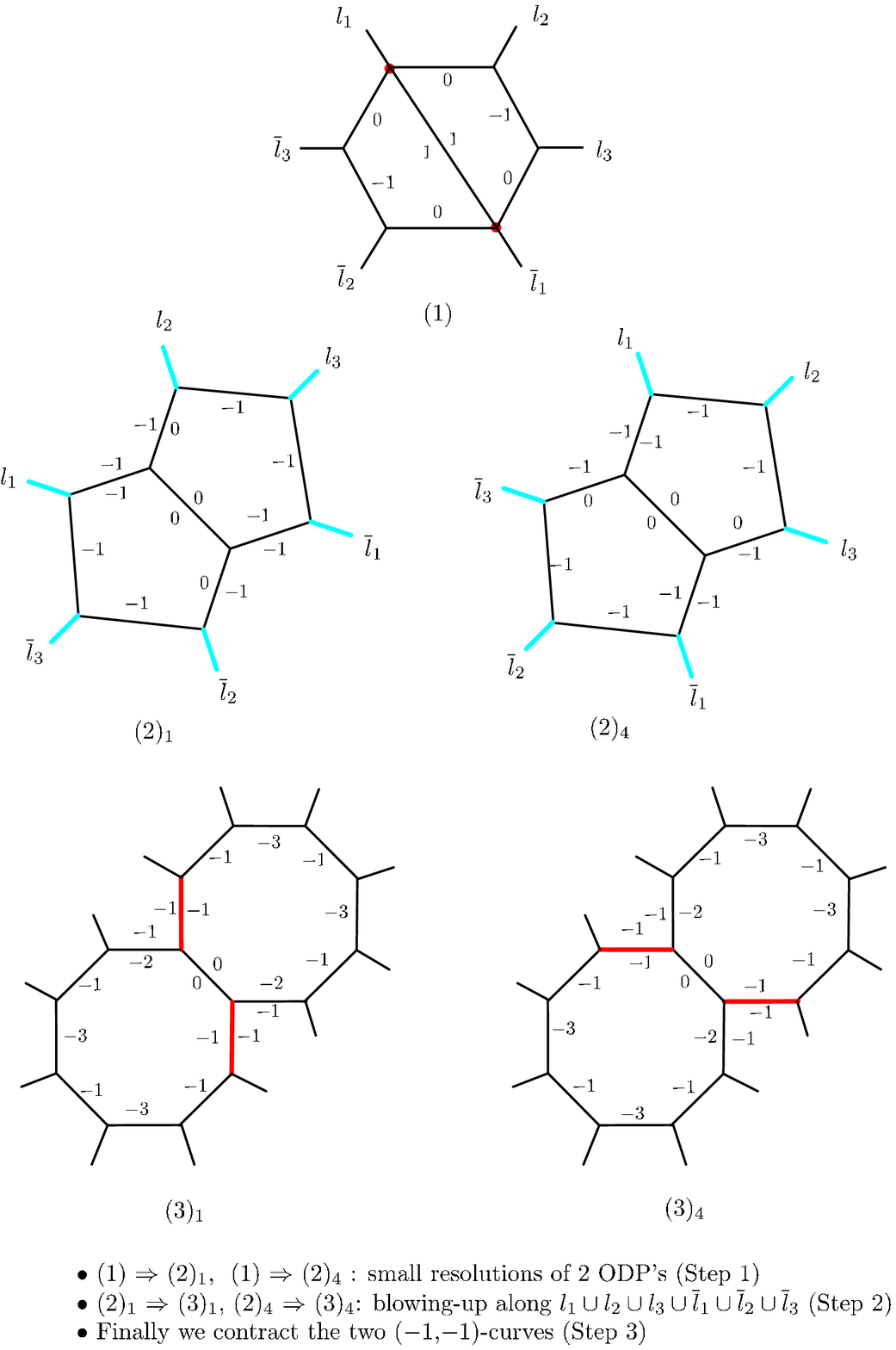}
\caption{The transformations for the reducible fiber over $\lambda_1$ and $\lambda_4$}
\label{fig-bim14}
\end{figure}

\begin{figure}
\includegraphics{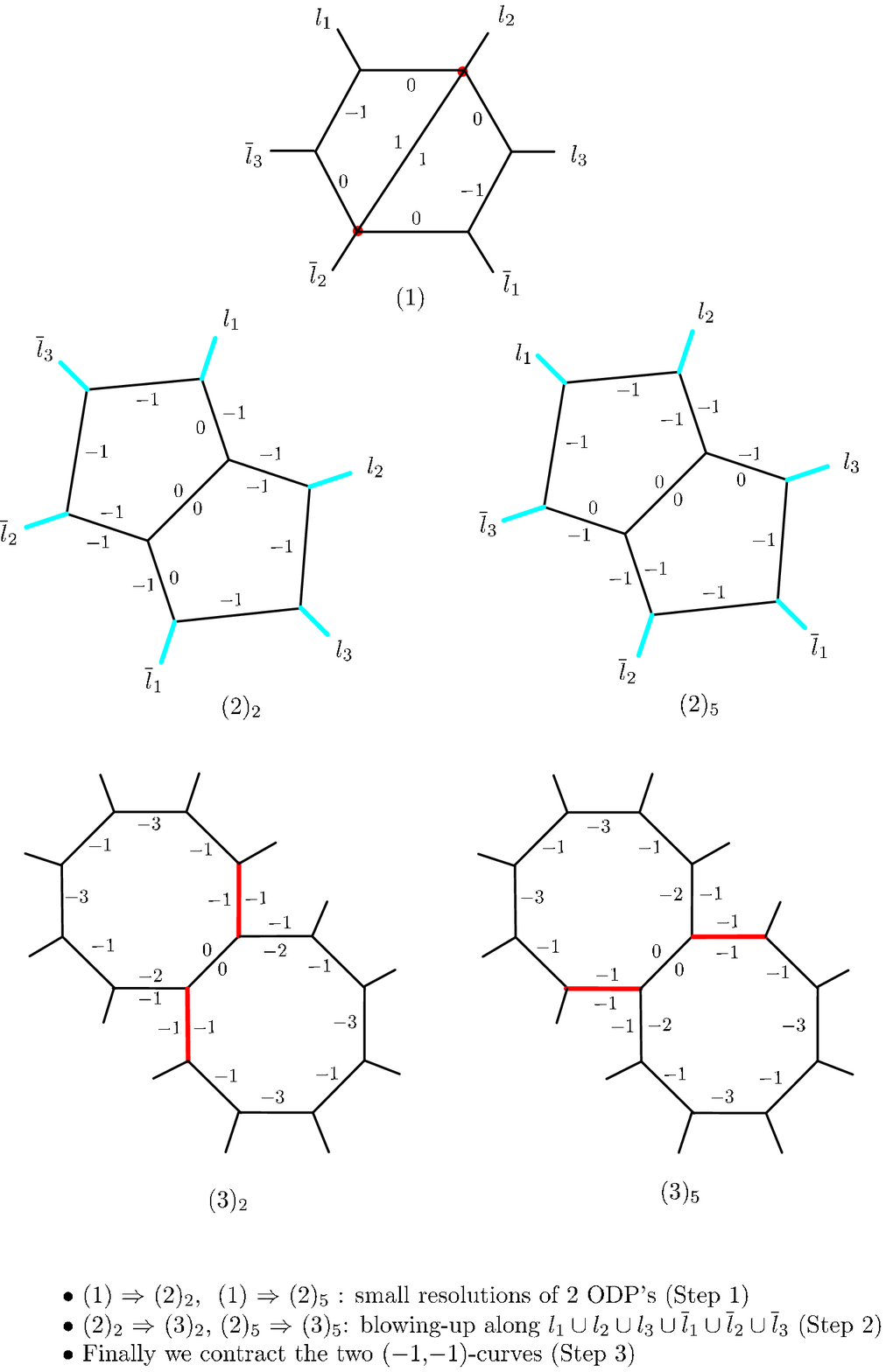}
\caption{The transformations for the reducible fiber over $\lambda_2$  and $\lambda_5$}
\label{fig-bim25}
\end{figure}

\begin{figure}
\includegraphics{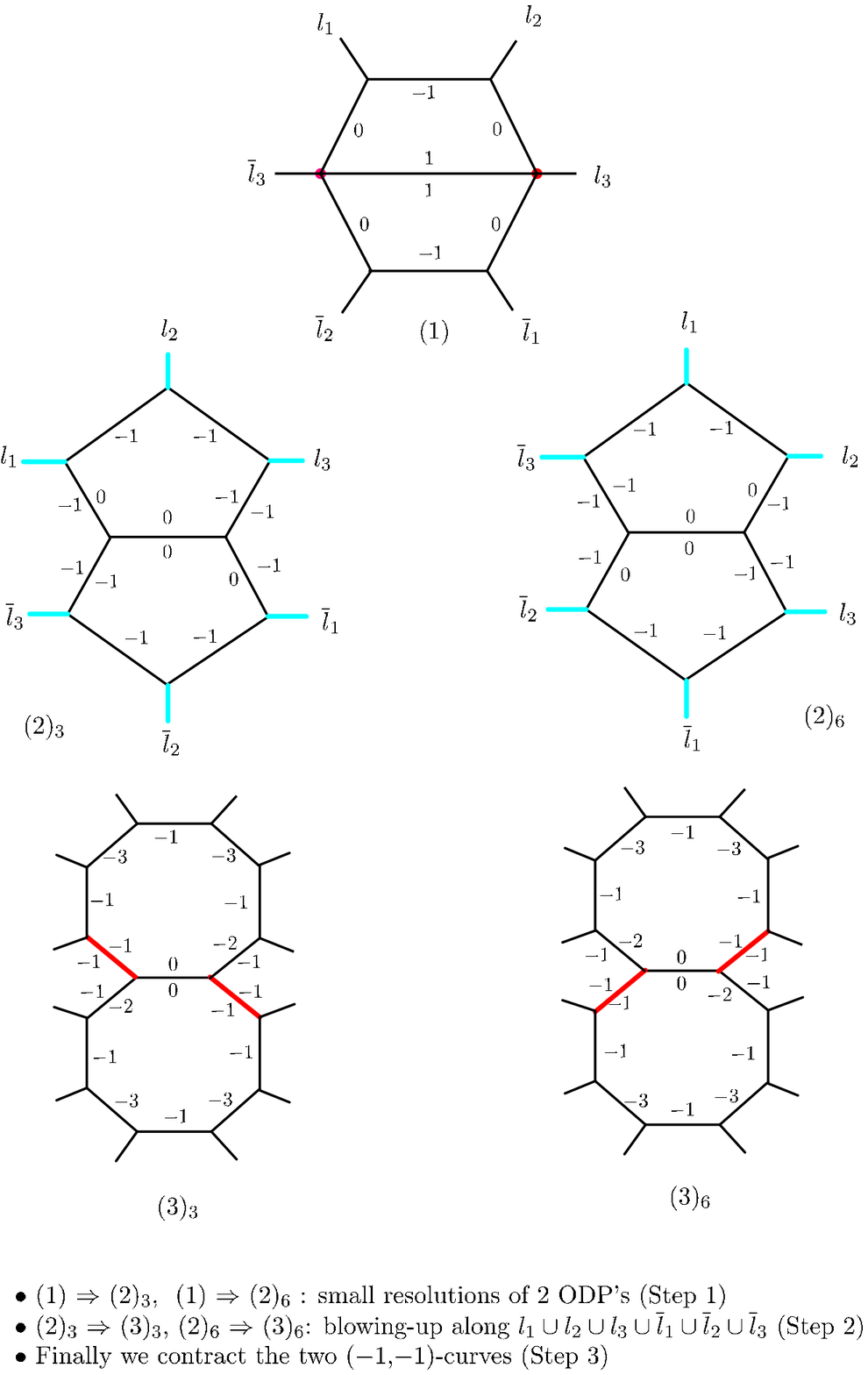}
\caption{The transformations for the reducible fiber over $\lambda_3$  and $\lambda_6$}
\label{fig-bim36}
\end{figure}

On the total space of the bundle $\mathbf P(\mathscr O(2)^{\oplus 6}\oplus \mathscr O)\ra\Lambda$ we consider the $G$-action 
$$
(x_0,x_1,x_2,x_3,x_4,x_5,x_6)\longmapsto (x_0, sx_1,s^{-1}x_2,tx_3,t^{-1}x_4,s^{-1}t^{-1}x_5,stx_6)\hspace{2mm}\text{for}\hspace{2mm}(s,t)\in G,
$$
where $x_i$, $1\le i\le 6$ is a fiber coordinate on the line bundle $\mathscr O(2)$ and $x_0$ is a fiber coordinate on $\mathscr O$.
We also consider  the real structure
$$
(x_0,x_1,x_2,x_3,x_4,x_5,x_6)\longmapsto
(\ol{x}_0,\ol{x}_2,\ol{x}_1,\ol{x}_4,\ol{x}_3,\ol{x}_6,\ol{x}_5).
$$
(These are the ones obtained in Lemma \ref{lemma-01} and its proof.)
Then the following statements are either obvious from the results in the previous section or verified by explicit computations:

\vsp
(a) If $\lambda\neq\lambda_i$ for $1\le i\le 6$, the fiber $p_1^{-1}(\lambda)$ is a smooth toric surface with $K^2=6$.
On each fiber, there exist precisely six $G$-fixed points and they form 6 global sections of $p_1:Z_1\ra\mathbf P^1$. 
Because these 6 sections consist of 3 conjugate pairs, we write $l_1,l_2,l_3,\ol{l}_1,\ol{l}_2,\ol{l}_3$ for these sections (Figures \ref{fig-bim14}--\ref{fig-bim36}, (1)).
The union of these 6 sections coincides with the fixed locus of the  $G$-action on $X_1$.

(b) For $1\le i\le 6$, $p^{-1}_1(\lambda_i)$ is a reducible toric surface consisting of 2 irreducible components. 
These 2 components form a conjugate pair and each components are smooth toric cubic surfaces contained  in some linear subspace $\mathbf P^4$ in $\mathbf P^6$.
The intersections of these cubic surfaces are conics in some linear subspace $\mathbf P^2\subset \mathbf P^6$.
(These 6 conics will become $G$-invariant twistor lines in $Z$.)

(c) The singular locus of $X_1$ consists of 12 isolated ordinary double points.
These consist of 6 conjugate pairs and each pair is precisely the pair of $G$-fixed points on the conics in (b).

\vsp
In the following we apply  simple birational transformations to $X_1$ in order to obtain the required twistor spaces.
All transformations are performed in such a way that the $G$-action and the real structure are preserved.
So in each steps below we do not mention about them.

\vsp\noindent
$\bullet$ {\bf Step 1} (Small resolutions of all ODP's in $X_1$).
For each ODP of $X_1$, there are 2 choices of small resolutions. 
They are uniquely determined once we specify which irreducible components are blown-up at the point.
So let $X_2\ra X_1$ be the small resolutions as illustrated in Figures  \ref{fig-bim14}--\ref{fig-bim36}, (1) $\Ra$ (2), and $p_2:X_2\ra \Lambda$  the natural projection.
$X_2$ is non-singular.
We use the same symbols $l_i$ and $\ol{l}_i$ ($1\le i\le 3)$ for the sections of $p_2$ coming from the 6 sections in $X_1$.

\vsp\noindent
$\bullet$ {\bf Step 2} (Blowing-ups along 6 sections)
Let $X_3\ra X_2$ be the blowing-up along $l_1\cup l_2\cup l_3\cup \ol{l}_1\cup \ol{l}_2\cup \ol{l}_3$ and $p_3:X_3\ra\Lambda$ the natural projection
(Figures  \ref{fig-bim14}--\ref{fig-bim36}, $(2)\Ra(3)$).
Obviously $p_3$ has 6 reducible fibers and on each of them there is a unique pair of $(-1,-1)$-curves contained in each reducible fiber (Figures  \ref{fig-bim14}--\ref{fig-bim36}, $(3)$).

\vsp\noindent
$\bullet$ {\bf Step 3}  (Contractions of redundant 12 curves into ODP's).
Let $X_3\ra \hat Z$ be the contraction of these twelve $(-1,-1)$-curves into ordinary double points, and $f:\hat Z\ra\Lambda$ be the natural projection.
Any fiber of $f$ has a cycle of $G$-invariant rational curves consisting of 12 irreducible components.

\vsp
Then we have the following 

\begin{thm}\label{thm-02}
Let $g$ be the Joyce metric on $4\mathbf P^2$ whose $K$-action is Type II and whose conformal invariant is $\{\lambda_1<\cdots<\lambda_6\}$.
Let $Z$ be the twistor space of \,$g$ and $X_1$ the anticanonical model of\, $Z_1$ defined by the equations \eqref{defeq02-1}--\eqref{defeq02-9}.
Let $\hat Z$ be the 3-fold constructed from $X_1$ by Step 1--Step 3 above.
Let $\tilde Z\ra Z$ be the blowing-up along the cycle $C$ which is the base curve of the system $|(-1/2)K_Z|$.
Then there exists an isomorphism between $\tilde Z$ and  $\hat Z$ which commutes with the real structure.

\end{thm}

\noindent
Proof.
By using our explicit constructions, it can be verified that all the fibers of the 2 fiber spaces $\hat Z\ra\mathbf P^1$ and $\tilde Z\ra\mathbf P^1$ are isomorphic as toric surfaces.
It can also be checked  that the real structure on  $\hat Z$ has no fixed points and exchanges $C_i$ and $\ol{C}_i$ for $1\le i\le 6$.
So by Fujiki's proof of \cite[Theorem 8.1]{F00}, to show the existence of an isomorphism between $\hat Z$ and $\tilde Z$ compatible with the real structures, it suffices to prove  that there exists a real smooth rational curves $\hat L$ in $\hat Z$ satisfying  (i) $\hat L$ is two to one over $\Lambda$ and the 2 branch points of $\hat L\ra\Lambda$ do not on the real locus on $\Lambda$, and (ii) $\hat L$ is disjoint from the $G$-invariant cycle of 12 rational curves in any fibers.
For this, it suffices to take a general twistor line $L\subset Z$ which is disjoint from the cycle $C$ and letting $\hat L$ to be the strict transform of the image $\Phi_1(L)$ into $\hat Z$.
In fact, (i) is satisfied by the same reasoning as that in the last paragraph of \cite[Theorem 3.2]{H-1}.
Moreover, since $\Phi_1$ gives a biholomorphic map between the unions of all 2-dimensional $G$-orbits in $Z$ and $X_1$, 
and since Step 1--3 do not touch the latter union of $G$-orbits in $X_1$, 
$\hat L$ is disjoint from the $G$-invariant cycles in any fiber.
Thus $\hat L$ satisfies (ii) as well.
\proofend

\begin{rmk}{\rm
The above proof is pretty easier than that of the corresponding result (Theorem 3.1) in \cite{H-1}.
There, we have to take account of  the effect of  `Blowing-down of the redundant divisors' (Step 5 in \cite[Section 3]{H-1}), to obtain the above property (ii) for $\hat L$.
Namely, there was a possibility that the property (ii) is lost under the blowing-downs of the divisors.
 In the present case (Type II) we do not blow-down any divisors through Steps 1--3 and this makes the proof much simpler.}\end{rmk}

\newpage

\end{document}